\documentclass{tac}

\usepackage{amssymb,amsmath}
\usepackage{tikz-cd}

\usepackage{xy}

\input diagxy

\usepackage[colorlinks=true]{hyperref}
\hypersetup{allcolors=[rgb]{0.1,0.1,0.4}}

\author{Daniel Graves}

\address{School of Mathematics and Statistics, University of Sheffield, Sheffield, S3 7RH, UK
}

\title{PROPs for involutive monoids and involutive bimonoids}

\copyrightyear{2020}


\eaddress{dan.graves92@gmail.com}

\renewcommand{\phi}{\varphi}
\newcommand{\cyc}{C_2}
\newcommand{\fas}{\mathcal{F}(as)}
\newcommand{\ifas}{\mathcal{IF}(as)}
\newcommand{\C}{\mathbf{C}}
\newcommand{\llb}{\left\lbrace}
\newcommand{\rrb}{\right\rbrace}
\renewcommand{\ul}{\underline}
\newcommand{\llangle}{\left\langle}
\newcommand{\rrangle}{\right\rangle}
\newcommand{\mr}{\mathrm}

\newtheorem{theorem}{Theorem}

\newtheoremrm{rem}{Remark}

\theoremstyle{definition}

\mathrmdef{Hom}
\mathbfdef{Set}

\begin{document}

\maketitle
\begin{abstract}
The category of involutive non-commutative sets encodes the structure of an involution compatible with a (co)associative (co)multiplication. We prove that the category of involutive bimonoids in a symmetric monoidal category is equivalent to the category of algebras over a PROP constructed from the category of involutive non-commutative sets.
\end{abstract}


\section*{Introduction}
\label{sec-Introduction}
The categorification of algebras over a unital commutative ring $k$ to algebras over a PROP was first introduced by Markl in order to study the deformation theory of algebras \cite{markl}. In that paper he defined PROPs, in terms of generators and relations, whose categories of algebras are equivalent to the category of associative algebras, the category of commutative algebras and the category of bialgebras over $k$ \cite[Examples 2.5, 2.6 and 2.7]{markl}. 

Pirashvili \cite{Pir-PROP} gave an explicit description of a PROP that categorified associative algebras, commutative algebras and bialgebras in the category of vector spaces over a field. This PROP is constructed from the category of non-commutative sets, introduced by Feigin and Tsygan \cite[A10]{FT}, using the generalized Quillen $Q$-construction of Fiedorowicz and Loday \cite[2.5]{FL}. An alternative approach, using distributive laws for PROPs, was given by Lack \cite[Section 5]{Lack}. In this setting, Pirashvili's PROP is described as a composite constructed from the PRO of finite ordinals and its opposite category and the result is shown to be more general, holding for bimonoids in a symmetric monoidal category.

In this paper we combine both of these methods. We introduce the PROP of \emph{involutive non-commutative sets}, denoted $\ifas$. Using the machinery of \cite{Lack} we describe a composite PROP constructed from $\ifas$ whose algebras in a symmetric monoidal category are the involutive bimonoids. 

The paper is organized as follows. In Section \ref{alg-sec} we recall the definition of involutive bimonoid in a symmetric monoidal category. In Section \ref{PROP-sec} we recall the definitions of PRO and PROP, together with some examples. In Section \ref{ifas-sec} we define the PROP of involutive non-commutative sets, $\ifas$, and prove that the category of algebras of $\ifas$ in a symmetric monoidal category $\C$ is equivalent to the category of involutive monoids in $\C$. In Section \ref{double-sec} we construct a double category from $\ifas$ whose bimorphisms encode the compatibility conditions for an involutive bimonoid. In Section \ref{bimon-sec} we construct a composite PROP, in the sense of \cite[Section 4]{Lack}, from the PROP $\ifas$ and its opposite category. This can be seen as a composite of the PROPs for involutive monoids and involutive comonoids described in Section \ref{ifas-sec}, where the compatibility of the two is encoded in the double category described in Section \ref{double-sec}. We prove that the category of algebras of this composite PROP in a symmetric monoidal category $\C$ is equivalent to the category of involutive bimonoids in $\C$.

\section*{Acknowledgements}
I would like to thank Sarah Whitehouse for all her guidance and encouragement. I am grateful to James Brotherston, James Cranch and the anonymous referees for their helpful suggestions. 

\section{Involutive monoids and involutive bimonoids}
\label{alg-sec}
\definition
A monoid $M$ in a symmetric monoidal category $\C$ is called \emph{involutive} if it comes equipped with a  monoid morphism $j\colon M\rightarrow M^{op}$ satisfying $j^2=id_M$. We denote the category of involutive monoids and involution-preserving morphisms by $\mathbf{IMon}\left(\C\right)$.

We denote the category of involutive comonoids in $\C$ by $\mathbf{IComon}\left(\C\right) = \mathbf{IMon}\left(\C^{op}\right)^{op}$.

A bimonoid $B$ in $\C$ is said to be involutive if it comes equipped with a bimonoid morphism $j\colon B\rightarrow B^{op, cop}$ such that $j^2=id_B$. We denote the category of involutive bimonoids in $\C$ by $\mathbf{IBimon}\left(\C\right)$.
\enddefinition

\section{PROs and PROPs}
\label{PROP-sec}

\definition
\label{sets-defn}
For $n\geqslant 1$ we define $\ul{n}$ to be the set $\llb 1,\dotsc ,n\rrb$. We define $\ul{0}=\emptyset$.
\enddefinition

\definition
A PRO $\mathbb{T}$ is a strict monoidal category whose objects are the sets $\ul{n}$ for $n\geqslant 0$ and whose tensor product is given by the disjoint union. For a monoidal category $\C$, an \emph{algebra of $\mathbb{T}$ in $\C$} is a strict monoidal functor $\mathbb{T}\rightarrow \C$.
\enddefinition

\definition
A \emph{PROP} $\mathbf{P}$ is a symmetric strict monoidal category whose objects are the sets $\ul{n}$ for $n\geqslant 0$ with tensor product given by the disjoint union. For a symmetric monoidal category $\mathbf{C}$, a \emph{$\mathbf{P}$-algebra in $\mathbf{C}$} is a symmetric strict monoidal functor $\mathbf{P}\rightarrow\mathbf{C}$. We denote the category of $\mathbf{P}$-algebras in $\mathbf{C}$ and natural transformations by $\mathbf{Alg}\left(\mathbf{P},\mathbf{C}\right)$.
\enddefinition

\example
We denote by $\mathbb{D}$ the PRO of finite ordinals and order-preserving maps as in \cite[2.2]{Lack}. For a strict monoidal category $\C$, an algebra of $\mathbb{D}$ in $\C$ is a monoid in $\C$, see \cite[VII 5]{ML}.
\endexample

\definition
Let $C_2=\left\langle t\mid t^2=1\right\rangle$. Let $\mathbb{C}_2$ be the PRO such that $\mathrm{Hom}_{\mathbb{C}_2}\left(\ul{n} , \ul{m}\right)$ is empty if $n\neq m$ and $\mathrm{Hom}_{\mathbb{C}_2}\left(\ul{n} , \ul{n}\right) = C_2^n$. The disjoint union of morphisms corresponds to the product of group elements.
\enddefinition

\example
Following \cite[2.4]{Lack}, let $\mathbb{P}$ denote the PRO of finite sets and bijections.
\endexample

\remark
\label{Prop-remark}
An equivalent definition of a PROP is as a PRO $\mathbb{T}$ with a map of PROs $\mathbb{P}\rightarrow \mathbb{T}$. Therefore $\mathbb{P}$ is a PROP. Given a PROP $\mathbf{P}$ we will denote its underlying PRO by $\mathbf{P}_0$. 
\endremark

\example
Following \cite[5.1, 5.2]{Lack}, we denote by $\mathbb{F}$ the PROP of finite sets and finite set maps. For a symmetric monoidal category $\C$, the category $\mathbf{Alg}\left(\mathbb{F} , \C\right)$ is equivalent to the category of commutative monoids in $\C$ and $\mathbf{Alg}\left(\mathbb{F}^{op} , \C\right)$ is equivalent to the category of cocommutative comonoids in $\C$.
\endexample

We can form new PROs and PROPs via the notion of a \emph{distributive law} as defined in \cite[Section 3]{Lack} and \cite[Section 2]{RW}. In particular we will form a composite PROP from $\mathbb{D}$, $\mathbb{P}$ and $\mathbb{C}_2$ whose structure is that of the hyperoctahedral category defined by Fiedorowicz and Loday \cite[Section 3]{FL}. We begin by recalling the definition of the hyperoctahedral groups.

\definition
\label{hyp-oct-defn}
For $n\geqslant 1$, the \emph{hyperoctahedral group} $H_{n}$ is defined to be the semi-direct product $\cyc^{n} \rtimes \Sigma_{n}$ where $\Sigma_{n}$ acts on $\cyc^{n}$ by permuting the factors.
\enddefinition

\example 
A pair $\left(x, f\right)$ where $f\in \mathrm{Hom}_{\mathbb{F}}\left(\ul{n} , \ul{m}\right)$ and $x\in \mathrm{Hom}_{\mathbb{C}_2}\left(\ul{m} , \ul{m}\right)$ determines a unique pair $\left(f , x^{\prime}\right)$ where $f$ has remained unchanged and, if $x=\left(g_1,\dotsc , g_m\right)$, $x^{\prime}=\left(g_{f(1)},\dotsc , g_{f(n)}\right)\in \mathrm{Hom}_{\mathbb{C}_2}\left(\ul{n} , \ul{n}\right)$. A straightforward check of the relations in \cite[2.4]{RW} and \cite[3.7]{Lack} shows that this defines a distributive law $\mathbb{C}_2 \otimes \mathbb{F}\rightarrow \mathbb{F}\otimes \mathbb{C}_2$ which is compatible with the monoidal structures of $\mathbb{F}$ and $\mathbb{C}_2$.

By \cite[Theorem 3.8]{Lack}, $\mathbb{F}\otimes \mathbb{C}_2$ is a PRO such that morphisms in $\mathrm{Hom}_{\mathbb{F}\otimes \mathbb{C}_2}\left(\ul{n} , \ul{m}\right)$ can be written uniquely as pairs $\left(f,x\right)$ with $x\in \mathrm{Hom}_{\mathbb{C}_2}\left(\ul{n} , \ul{n}\right)$ and $f\in \mathrm{Hom}_{\mathbb{F}}\left(\ul{n} , \ul{m}\right)$ with composition determined by the distributive law. In fact, $\mathbb{F}\otimes \mathbb{C}_2$ has a canonical PROP structure induced from the PROP structure on $\mathbb{F}$.
\endexample

\example
A pair $\left(x , \sigma\right)$ where $\sigma\in \mathrm{Hom}_{\mathbb{P}}\left(\ul{n} ,\ul{n}\right)$ and $x\in \mathrm{Hom}_{\mathbb{C}_2}\left(\ul{n} , \ul{n}\right)$ determines a unique pair $\left(\sigma , x^{\prime}\right)$ where $\sigma$ has remained unchanged and, if $x=\left(g_1, \dotsc , g_n\right)$, $x^{\prime}=\left(g_{\sigma(1)},\dotsc , g_{\sigma(n)}\right)$.

Similarly to the previous example, this is a distributive law compatible with the monoidal structure of $\mathbb{P}$ and $\mathbb{C}_2$ and we have a PRO $\mathbb{P}\otimes \mathbb{C}_2$.

We observe that $\mathrm{Hom}_{\mathbb{P}\otimes \mathbb{C}_2}\left(\ul{n} , \ul{n}\right)$ is isomorphic to $H_n$, the $n^{th}$ hyperoctahedral group as defined in Definition \ref{hyp-oct-defn}. We therefore write $\mathbb{H}=\mathbb{P}\otimes \mathbb{C}_2$. We note that $\mathbb{H}$ has a canonical PROP structure induced from $\mathbb{P}$. We refer to $\mathbb{H}$ as the PROP of hyperoctahedral groups.
\endexample

\example
Given a pair $\left(g, \phi\right)$ where $\phi \in \mathrm{Hom}_{\mathbb{D}}\left(\ul{n} , \ul{m}\right)$ and $g\in \mathrm{Hom}_{\mathbb{H}}\left(\ul{m} , \ul{m}\right)$ there is a unique pair $\left(g_{\star}\left(\phi\right),\phi^{\star}(g)\right)$ where $\phi^{\star}(g)\in \mathrm{Hom}_{\mathbb{H}}\left(\ul{n} ,\ul{n}\right)$ and $g_{\star}(\phi)\in \mathrm{Hom}_{\mathbb{D}}\left(\ul{n} , \ul{m}\right)$, as constructed in \cite[3.1]{FL}. The fact that these assignments satisfy the relations of \cite[2.4]{RW} follows from the fact that they satisfy the relations of a crossed simplicial group given in \cite[1.6]{FL} and a routine check shows that they respect the monoidal structures of $\mathbb{D}$ and $\mathbb{H}$.

We therefore have a PRO $\mathbb{D}\otimes \mathbb{H}$ defined similarly to the examples above. In fact, $\mathbb{D}\otimes \mathbb{H}$ has a canonical PROP structure induced from $\mathbb{H}$.
\endexample

\proposition
\label{comp-prop}
Let $\C$ be a symmetric monoidal category. There is an equivalence of categories
\[\mathbf{Alg}\left( \mathbb{D}\otimes \mathbb{H} , \C\right) \simeq \mathbf{IMon}\left(\C\right).\]
\endproposition
\proof
By \cite[3.10]{Lack}, a $\left(\mathbb{D}\otimes \mathbb{H}\right)_0$-algebra structure on an object $M$ of $\C$ consists of a $\mathbb{D}$-algebra structure and a $\mathbb{H}$-algebra structure subject to a compatibility condition. A $\mathbb{D}$-algebra structure is a monoid structure. A $\mathbb{H}$-algebra is an object $M$ together with a morphism $j\colon M\rightarrow M$ satisfying $j^2=id_M$ and, for each element $g\in H_n$, a morphism $M^{\otimes n}\rightarrow M^{\otimes n}$ given by applying $j$ to the tensor factors according to the element of $C_2^n$ followed by an isomorphism determined by the element of $\Sigma_n$. Arguing analogously to \cite[5.5]{Lack} a $\left(\mathbb{D}\otimes \mathbb{H}\right)_0$-algebra structure is a $\mathbb{D}\otimes \mathbb{H}$-algebra structure if and only if the only isomorphisms $M^{\otimes n}\rightarrow M^{\otimes n}$ are those induced from the symmetry isomorphisms. The compatibility condition is precisely the condition requiring $j$ to be an involution compatible with the monoid structure. Finally, a morphism in $\C$ is a map of involutive monoids if and only if it respects the $\mathbb{D}$-algebra structure and the $\mathbb{H}$-algebra structure. By \cite[3.12]{Lack}, this is true if and only if it respects the $\mathbb{D}\otimes \mathbb{H}$-algebra structure.
\endproof

\remark
This result tells us that the PROP governing the structure of an involutive monoid can be thought of as a composite of PROPs governing the structure of a monoid and the structure of an involution respectively. We will give an explicit description of this category, where the technicalities of distributive laws are distilled into data on the preimages of set maps in Section \ref{ifas-sec}.

We have chosen to emphasize the connection between involutive monoids and the category associated to the hyperoctahedral crossed simplicial group. One advantage of this approach is that the distributive laws employed are already well-known, being the composition in hyperoctahedral groups and the hyperoctahedral category. It is also an interesting new application of the hyperoctahedral crossed simplicial group: the other known applications are found in the field of equivariant stable homotopy theory!

An alternative method of proof would be to begin with the composite PROP $\mathbb{D}\otimes \mathbb{P}$ of \cite[3.14]{Lack}, define a distributive law between this and the PRO $\mathbb{C}_2$ and to analyse the resulting composite.
\endremark

\corollary
\label{comp-cor}
There is an equivalence of categories
\[\mathbf{Alg}\left( \left(\mathbb{D}\otimes \mathbb{H}\right)^{op} , \C\right) \simeq \mathbf{IComon}\left(\C\right).\]
\endcorollary
\proof
We observe that
\[\mathbf{IComon}\left(\C\right)=\mathbf{IMon}\left(\C^{op}\right)^{op}\simeq \mathbf{Alg}\left( \left(\mathbb{D}\otimes \mathbb{H}\right) , \C^{op}\right)^{op}=\mathbf{Alg}\left( \left(\mathbb{D}\otimes \mathbb{H}\right)^{op} , \C\right)\]
as required.
\endproof

\proposition
\label{comm-mon-prop}
Let $\C$ be a symmetric monoidal category. There is an equivalence of categories between $\mathbf{Alg}\left(\mathbb{F}\otimes \mathbb{C}_2 , \C\right)$ and the category of involutive commutative monoids in $\C$. The category $\mathbf{Alg}\left(\left(\mathbb{F}\otimes \mathbb{C}_2\right)^{op} , \C\right)$ is equivalent to the category of involutive cocommutative comonoids in $\C$.
\endproposition
\proof
The proof of the first equivalence is similar to Proposition \ref{comp-prop}. We note that an $\mathbb{F}$-algebra structure is a commutative monoid structure. A $\mathbb{C}_2$-algebra structure consists of an object $M$ in $\C$ together with a morphism $j\colon M\rightarrow M$ satisfying $j^2=id_M$ and for each element of $C_2^n$ a morphism $M^{\otimes n}\rightarrow M^{\otimes n}$ defined by applying $j$ to the tensor factors according to the element of $C_2^n$. An $\left(\mathbb{F}\otimes \mathbb{C}_2\right)_0$-structure is a $ \mathbb{F}\otimes \mathbb{C}_2$-structure if and only if the only isomorphisms $M^{\otimes n}\rightarrow M^{\otimes n}$ are those induced from the symmetry isomorphisms. The compatibility condition in this case is the condition that requires $j$ to be an involution compatible with a commutative monoid structure. Finally we note that a morphism in $\C$ is a map of involutive commutative monoids if and only if it preserves both the $\mathbb{F}$-algebra structure and the $\mathbb{C}_2$-algebra structure. The second equivalence follows a similar argument to Corollary \ref{comp-cor}.
\endproof

\section{The PROPs $\ifas$ and $\mathcal{IF}$}
\label{ifas-sec}
In the previous section we described a PROP for involutive monoids as a composite. In this section we provide an explicit description of this PROP, called the category of \emph{involutive non-commutative sets}. A variant of this category first appeared in the author's thesis \cite[Part V]{DMG}. This category takes the technicalities of the composition of pairs defined via a distributive law and presents it as simple structure on the preimages of maps of finite sets. We shall also see, in Section \ref{double-sec}, that we can construct a double category from the category of involutive non-commutative sets whose bimorphisms encode the structure of an involutive bimonoid.

The PROP of involutive, non-commutative sets, $\ifas$ will have as objects the sets $\ul{n}$ of Definition \ref{sets-defn} for $n\geqslant 0$. An element $f\in \mr{Hom}_{\ifas}\left(\ul{n},\ul{m}\right)$ will be a map of sets such that the preimage of each singleton $i\in \ul{m}$ is a totally ordered set such that each element comes adorned with a superscript label from the group $\cyc=\llangle t\mid t^2=1\rrangle$. Note that for $m\geqslant 1$, the set $\mr{Hom}_{\ifas}\left(\ul{0},\ul{m}\right)$ will be the singleton set consisting of the unique set map $\emptyset \rightarrow \ul{m}$ and $\mr{Hom}_{\ifas}\left(\ul{m},\ul{0}\right)$ will be the empty set. 

\remark
Henceforth we will say that a morphism in $\ifas$ is a map of sets together with a \emph{labelled, ordered set} for each preimage. In particular, note that we will use \emph{preimage} to mean preimage of a singleton.
\endremark

\example
\label{ifas-morph-eg}
Let $f\in\mr{Hom}_{\ifas}\left(\ul{5},\ul{4}\right)$ have underlying map of sets

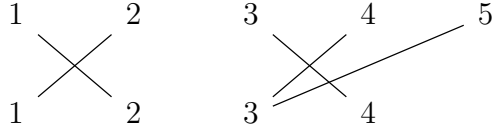
\begin{center}
\begin{tikzcd}
1\arrow[dr,dash] & 2\arrow[dl,dash] &3\arrow[dr,dash]&4\arrow[dl,dash]&5\arrow[dll,dash]\\
1&2&3&4
\end{tikzcd}
\end{center}
with the following labelled, ordered sets as preimages:
\[f^{-1}(1)=\llb 2^1\rrb, \quad f^{-1}(2)=\llb 1^t\rrb,\quad f^{-1}(3)=\llb 4^t<5^1\rrb \text{ and } f^{-1}(4)=\llb 3^t\rrb.\] 
\endexample

We will denote composition in $\ifas$ by $\bullet$ in order to distinguish from the composition of maps of sets. In particular, we use $\circ$ for two morphisms in $\ifas$ if we are referring to the composite of the underlying maps of sets. In order to ease notation we have chosen not to introduce notation for the forgetful functor $\ifas \rightarrow \mathbf{Set}$.

Let $f_1\in \mr{Hom}_{\ifas}\left(\ul{n},\ul{m}\right)$ and $f_2\in \mr{Hom}_{\ifas}\left(\ul{m},\ul{l}\right)$. In order to define the composite $f_2\bullet f_1\in \mr{Hom}_{\ifas}\left(\ul{n},\ul{l}\right)$ we must provide a map of sets and describe the labelled total orderings on each of the preimages.

As a map of sets, $f_2\bullet f_1$ is the composite of the underlying map of sets $f_2\circ f_1$. In order to specify a labelled, ordered set for the preimage of each singleton in $\ul{l}$ under the composite we first make a definition.

\definition
We define an action of $\cyc$, which will be denoted by a superscript, on finite, ordered sets with $\cyc$-labels by
\[ \llb j_1^{\alpha_{j_1}}<\cdots <j_r^{\alpha_{j_r}}\rrb^t=\llb j_r^{t\alpha_{j_r}}<\cdots <j_1^{t\alpha_{j_1}}\rrb.\]
That is, we invert the ordering and multiply each label by $t\in \cyc$.
\enddefinition

\definition
\label{ifas-comp-defn}
Let $f_1\in \mr{Hom}_{\ifas}\left(\ul{n},\ul{m}\right)$ and $f_2\in \mr{Hom}_{\ifas}\left(\ul{m},\ul{l}\right)$. We define $f_2\bullet f_1\in \mr{Hom}_{\ifas}\left(\ul{n},\ul{l}\right)$ to have underlying map of sets $f_2\circ f_1$.
We define the labelled totally ordered set $(f_2\bullet f_1)^{-1}(i)$ to be the ordered disjoint union of labelled, ordered sets
\[ \coprod_{j^{\alpha_j}\in f_2^{-1}(i)} f_1^{-1}\left(j\right)^{\alpha_j}.\]
\enddefinition

\definition
The \emph{PROP of involutive, non-commutative sets}, $\ifas$, has as objects the sets $\ul{n}$ of Definition \ref{sets-defn} for $n\geqslant 0$. An element of $\mr{Hom}_{\ifas}(\ul{n},\ul{m})$ is a map of sets with a total ordering on each preimage such that each element of the domain comes adorned with a superscript label from the group $\cyc$. Composition of morphisms is as defined in Definition \ref{ifas-comp-defn}. The symmetry isomorphisms are given by block permutations.
\enddefinition

\remark
For each $m\geqslant 1$, the set $\mr{Hom}_{\ifas}\left(\ul{0},\ul{m}\right)$ is the singleton set consisting of the unique set map $\emptyset \rightarrow \ul{m}$ and $\mr{Hom}_{\ifas}\left(\ul{m},\ul{0}\right)$ is the empty set.
\endremark

\remark
Recall the PROP of non-commutative sets, $\fas$, from \cite[Section 3]{Pir-PROP}. That is, the category whose objects are the sets $\ul{n}$ for $n\geqslant 0$ and whose morphisms are maps of sets with a total ordering on the preimage of each singleton in the codomain. We observe that $\fas$ is isomorphic to the subcategory of $\ifas$ which contains only the morphisms for which every label is $1\in \cyc$.
\endremark

\definition
\label{fund-morph-defn}
We define the \emph{fundamental morphisms} $m$, $u$ and $i$ of $\ifas$ as follows.
\begin{itemize}
\item Let $m\in \mr{Hom}_{\ifas}\left(\ul{2},\ul{1}\right)$ be defined by $m^{-1}(1)=\llb 1^1<2^1\rrb$,
\item let $u$ be the unique morphism in $\mr{Hom}_{\ifas}\left(\ul{0},\ul{1}\right)$ and
\item let $i\in \mr{Hom}_{\ifas}\left(\ul{1},\ul{1}\right)$ be defined by $i^{-1}(1)=\llb 1^t\rrb$.
\end{itemize}
\enddefinition

\remark
The morphism $m$ will encode the multiplication and comultiplication in a bimonoid, the morphism $u$ will encode the unit and counit and the morphism $i$ will encode the involution.
\endremark

\remark
We note that $\ifas$ contains the morphisms of the PRO of finite ordinals $\mathbb{D}$. These are the order-preserving maps of sets with the canonical total ordering on each preimage with each label being $1\in C_2$. Furthermore, $\ifas$ contains the morphisms of the PROP of hyperoctahedral groups $\mathbb{H}$. These are the bijections in $\ifas$.
\endremark

\proposition
\label{iso-prop}
There is an isomorphism of PROPs $\ifas\cong \mathbb{D}\otimes \mathbb{H}$. 
\endproposition
\proof
Consider the data of a morphism $f\in \mathrm{Hom}_{\ifas}\left(\ul{n} , \ul{m}\right)$. The total ordering data on preimages determines a unique bijection of the set $\ul{n}$ with a $C_2$-label for each singleton preimage determined by the labelling data. That is, the preimage data determines a unique bijection $g \in \mathrm{Hom}_{\ifas}\left(\ul{n} , \ul{n}\right)$.  There is then a unique morphism $\phi \in \mathrm{Hom}_{\ifas}\left(\ul{n} , \ul{m}\right)$ such that $\phi$ is order-preserving, with the canonical total ordering on each preimage, every label is $1\in C_2$ and $f=\phi \bullet g$. In other words, any morphism in $\ifas$ can be written uniquely as a composite of a morphism in $\mathbb{H}$ followed by a morphism in $\mathbb{D}$. 

For composable morphisms $f_1$ and $f_2$ in $\ifas$ write $f_2\bullet f_1= \phi_{2}\bullet g_{2}\bullet \phi_{1} \bullet g_{1}$. A straightforward check shows that the composite $g_2\bullet \phi_1$ in $\ifas$ is equal to the composite $g_{2\star}\left(\phi_1\right) \bullet \phi_1^{\star}\left(g_2\right)$, where $g_{2\star}\left(\phi_1\right) \in \mathrm{Hom}_{\ifas}\left(\ul{n} , \ul{m}\right)$ is an order-preserving map and $\phi_1^{\star}\left(g_2\right)\in \mathrm{Hom}_{\ifas}\left(\ul{n} ,\ul{n}\right)$ is a bijection, both maps being determined using the structure of the hyperoctahedral crossed simplicial group described in \cite[3.1]{FL}. It follows that there is an isomorphism of PROPs $\ifas \cong \mathbb{D} \otimes \mathbb{H}$ given by sending a morphism $f=\phi\bullet g$ in $\ifas$ to the pair $\left(\phi , g\right)$.
\endproof

\corollary
Let $\C$ be a symmetric monoidal category. There are equivalences of categories 
\[\mathbf{Alg}\left(\ifas , \C\right)\simeq \mathbf{IMon}\left(\C\right) \quad \text{and} \quad \mathbf{Alg}\left(\ifas^{op} , \C\right)\simeq \mathbf{IComon}\left(\C\right).\]
\endcorollary
\proof
This follows from Proposition \ref{iso-prop}, Proposition \ref{comp-prop} and Corollary \ref{comp-cor}.
\endproof

\definition
Let $\mathcal{IF}$ be the category whose objects are the sets $\ul{n}$ of Definition \ref{sets-defn} for $n\geqslant 0$. A morphism in $\mathcal{IF}$ is a map of sets such that the elements of the preimage of each singleton in the codomain come adorned with a label from $\cyc$. Composition is given by composition of set maps and multiplication of labels.
\enddefinition

\proposition
\label{if-prop}
There is an isomorphism of PROPs $\mathcal{IF}\cong \mathbb{F}\otimes \mathbb{C}_2$. 
\endproposition
\proof
The method of proof is similar to Proposition \ref{iso-prop}.
\endproof

\corollary
Let $\C$ be a symmetric monoidal category. The categories $\mathbf{Alg}\left(\mathcal{IF} , \C\right)$ and $\mathbf{Alg}\left(\mathcal{IF}^{op} , \C\right)$ are equivalent to the category of involutive monoids in $\C$ and the category of involutive cocommutative comonoids in $\C$ respectively.
\endcorollary
\proof
This follows from Proposition \ref{comm-mon-prop} and Proposition \ref{if-prop}.
\endproof

\section{Double categories}
\label{double-sec}
We construct a double category from $\ifas$. The bimorphisms of this double category precisely encode the structure of an involutive bimonoid in a symmetric monoidal category. This double category also possesses extra structure: it satisfies the star condition of \cite[2.3]{FL}. This extra structure will be used in Section \ref{bimon-sec} to construct a PROP which governs the structure of an involutive bimonoid. We also construct a double category from $\mathcal{IF}$ and two double categories that combine the structure of $\ifas$ and $\mathcal{IF}$ which will encode commutativity and cocommutativity.

Recall from \cite[Section 2.1]{FL} that a small \emph{double category} $\mathbf{D}$ consists of a set of objects, a set of horizontal morphisms, a set of vertical morphisms and a set of bimorphisms subject to natural composition identities.

\definition
The double category $\ifas_2$ has as objects the objects of $\ifas$. Furthermore, the sets of horizontal and vertical morphisms in $\ifas_2$ are both equal to the set of all morphisms in $\ifas$. A bimorphism in $\ifas_2$ is a not necessarily commutative square
\begin{center}
\begin{tikzcd}
\ul{n}\arrow[d, "\phi_1 ", swap]\arrow[r, "f_1"] & \ul{p}\arrow[d, "\phi "]\\
\ul{m}\arrow[r, "f", swap] & \ul{q}
\end{tikzcd}
\end{center}
of morphisms in $\ifas$ such that
\begin{itemize}
\item the underlying diagram of finite sets is a pullback square,
\item for all $x\in \ul{m}$ the map $\phi_1^{-1}(x)\rightarrow \phi^{-1}(f(x))$ induced by $f_1$ is an isomorphism of labelled, ordered sets and
\item for all $y\in \ul{p}$ the map $f_1^{-1}(y)\rightarrow f^{-1}(\phi(y))$ induced by $\phi_1$ is an isomorphism of labelled, ordered sets.
\end{itemize}
\enddefinition

\remark
The composition laws of a double category can be verified using the fact that the composite of pullback squares is itself a pullback square and using the composition rule for morphisms in $\ifas$ described in Definition \ref{ifas-comp-defn}.
\endremark

\definition
\label{fund-bimorph}
Let $B_1$, $B_2$, $B_3$, $B_4$ and $J$ denote the bimorphisms 
\begin{center}
\begin{tikzcd}
\ul{4}\arrow[d, swap, "{m^{\amalg 2} \bullet \tau_{2,3}}"]\arrow[r, "{m^{\amalg 2}}"] & \ul{2}\arrow[d, "m"]  & \ul{0}\arrow[d, swap, "{id_{\ul{0}}}"]\arrow[r, "{id_{\ul{0}}}"] & \ul{0}\arrow[d, "u"]  &  \ul{0}\arrow[d,swap, "{u^{2}}"]\arrow[r, "{id_{\ul{0}}}"] & \ul{0}\arrow[d, "u"] & \ul{0}\arrow[d, swap, "{id_{\ul{0}}}"]\arrow[r, "{u^{2}}"] &\ul{2}\arrow[d, "m"] &\ul{1}\arrow[d, swap, "{id_{\ul{1}}}"]\arrow[r, "{id_{\ul{1}}}"]& \ul{1}\arrow[d, "i"] \\
\ul{2}\arrow[r, "m", swap] & \ul{1} & \ul{0}\arrow[r, swap, "u"] & \ul{1}  & \ul{2}\arrow[r, "m",swap] & \ul{1} & \ul{0}\arrow[r, "u",swap] & \ul{1}& \ul{1}\arrow[r, swap, "i"] & \ul{1}
\end{tikzcd}
\end{center}
respectively in $\ifas_2$. Here $\tau_{2,3}$ is the transposition $(2\,\, 3)$ with the label $1\in \cyc$ for each preimage. We call $B_1$, $B_2$, $B_3$, $B_4$ and $J$ the \emph{fundamental bimorphisms of $\ifas_2$}.  
\enddefinition

\remark
The fundamental bimorphisms encode the compatibility conditions of an involutive bimonoid. The notation is chosen such that the bimorphisms $B_1$ to $B_4$ encode the compatibility conditions of a bimonoid and $J$ encodes the compatibility condition of an involution.
\endremark

\definition
The double category $\mathcal{IF}_2$ is defined similarly to $\ifas_2$; the objects are those of $\mathcal{IF}$, the sets of horizontal and vertical morphisms are the set of morphisms in $\mathcal{IF}$ and the bimorphisms are defined similarly to the bimorphisms of $\ifas_2$. 
\enddefinition

\definition
The double category $\mathcal{V}$ has as objects the objects of $\ifas$. The set of vertical morphisms is the set of morphisms in $\ifas$. The set of horizontal morphisms is the set of morphisms in $\mathcal{IF}$. The bimorphisms are defined similarly to those of $\ifas_2$ except that the horizontal morphisms are now in $\mathcal{IF}$.

The double category $\mathcal{H}$ is defined similarly; the set of horizontal morphisms is the set of morphisms in $\ifas$, the set of vertical morphisms is the set of morphisms in $\mathcal{IF}$ and the bimorphisms are defined similarly to those of $\ifas_2$ except that the vertical morphisms are in $\mathcal{IF}$.
\enddefinition

\remark
\label{uniq-bimorph-rem}
Recall from \cite[2.3]{FL} that a double category $\mathbf{D}$ is said to satisfy the \emph{star condition} if a horizontal morphism and a vertical morphism with the same codomain determine a unique bimorphism in $\mathbf{D}$.

Let $\mathbf{D}= \ifas_2$, $\mathcal{IF}_2$, $\mathcal{V}$ or $\mathcal{H}$. Given a horizontal morphism $f\colon \ul{m}\rightarrow \ul{q}$ and a vertical morphism $\phi \colon \ul{p}\rightarrow \ul{q}$ in $\mathbf{D}$ we determine a unique bimorphism by first taking the pullback of the underlying maps of sets. The resulting maps have a unique lift to the category $\mathbf{D}$ where the preimage data is induced from $f$ and $\phi$ using the conditions on bimorphisms. Therefore these four double categories satisfy the star condition.
\endremark

\section{Involutive bimonoids}
\label{bimon-sec}
We construct composite PROPs, in the sense of \cite[Section 4]{Lack}, from the PROPs $\ifas$, $\mathcal{IF}$ and their opposites. We prove that the category of algebras of the composite PROP constructed from $\ifas$ and its opposite in a symmetric monoidal category $\C$ is equivalent to the category of involutive bimonoids in $\C$.

\proposition
There exist composite PROPs 
\[\mathcal{IF}\otimes_{\mathbb{P}}\mathcal{IF}^{op}, \quad \mathcal{IF}\otimes_{\mathbb{P}}\ifas^{op}, \quad \ifas \otimes_{\mathbb{P}} \mathcal{IF}^{op}\quad \text{and} \quad \ifas \otimes_{\mathbb{P}} \ifas^{op}.\]
\endproposition
\proof
We provide the details for the case $\mathcal{IF}\otimes_{\mathbb{P}}\mathcal{IF}^{op}$. The others are similar. 

A pair $(\phi , f)$ where $f\in \mathrm{Hom}_{\mathcal{IF}}\left(\ul{m} ,\ul{q}\right)$ and $\phi \in \mathrm{Hom}_{\mathcal{IF}^{op}}\left(\ul{q} , \ul{p}\right)$ can be written as a diagram
\begin{center}
\begin{tikzcd}
& \ul{p}\arrow[d, "\phi "]\\
\ul{m}\arrow[r , "f" , swap] & \ul{q}
\end{tikzcd}
\end{center}
in $\mathcal{IF}$. By the star condition for the double category $\mathcal{IF}_2$, there exist unique morphisms $\phi_1$ and $f_1$ in $\mathcal{IF}$ forming a bimorphism
\begin{center}
\begin{tikzcd}
\ul{n}\arrow[d, "\phi_1 ", swap]\arrow[r, "f_1"] & \ul{p}\arrow[d, "\phi "]\\
\ul{m}\arrow[r, "f", swap] & \ul{q}
\end{tikzcd}
\end{center}
in $\mathcal{IF}_2$.

We observe that the assignment $\left(\phi , f\right)\mapsto \left(f_1 , \phi_1\right)$ defines a distributive law of PROs $ \mathcal{IF}^{op}\otimes \mathcal{IF}\rightarrow \mathcal{IF} \otimes \mathcal{IF}^{op}$, in the sense of \cite[3.6]{Lack}. The star condition, together with the composition rule for bimorphisms, ensures that the equations for a distributive law are satisfied and compatibility with the monoidal structure follows from the compatibility of the star condition with the disjoint union.

Furthermore, since both $\mathcal{IF}$ and $\mathcal{IF}^{op}$ are PROPs, we observe that $\mathcal{IF} \otimes \mathcal{IF}^{op}$ has a PROP structure.

A morphism in $\mathcal{IF} \otimes \mathcal{IF}^{op}$ from $\ul{n}$ to $\ul{m}$ can be written as a span
\begin{center}
\begin{tikzcd}
\ul{n} & \ul{p}\arrow[l, "\phi ", swap]\arrow[r, "f"] & \ul{m}.
\end{tikzcd}
\end{center}

Two spans
\begin{center}
\begin{tikzcd}
\ul{n} & \ul{p}\arrow[l, "\phi ", swap]\arrow[r, "f"] & \ul{m} & \text{and} & \ul{n} & \ul{p}\arrow[l, "\phi_1 ", swap]\arrow[r, "f_1"] & \ul{m}
\end{tikzcd}
\end{center}
are said to be \emph{equivalent} if there exists a bijection $h\colon \ul{p}\rightarrow \ul{p}$ in $\mathcal{IF}$ such that $\phi_1\circ h=\phi$ and $f_1\circ h=f$.

It follows from \cite[Theorem 4.6]{Lack} that $\mathcal{IF}\otimes_{\mathbb{P}}\mathcal{IF}^{op}$, that is the category obtained from $\mathcal{IF} \otimes \mathcal{IF}^{op}$ by identifying equivalent spans, is a composite PROP defined via a distributive law induced from the one defined for PROs.

The remaining three cases are similar, making use of the star condition from the double categories $\mathcal{V}$, $\mathcal{H}$ and $\ifas_2$ respectively.
\endproof

\definition
For ease of notation, let $Q=\ifas \otimes_{\mathbb{P}} \ifas^{op}$. Let $Q_{\mathcal{V}}=\mathcal{IF}\otimes_{\mathbb{P}}\ifas^{op}$. Let $Q_{\mathcal{H}}=\mathcal{IF}^{op}\otimes_{\mathbb{P}}\ifas$. Let $Q_{\mathcal{IF}}= \mathcal{IF}\otimes_{\mathbb{P}}\mathcal{IF}^{op}$.
\enddefinition

\theorem
\label{bimon-thm}
Let $\C$ be a symmetric monoidal category. There is an equivalence of categories
\[\mathbf{Alg}\left(Q , \C\right) \simeq \mathbf{IBimon}\left(\C\right).\]
\endtheorem
\proof
By \cite[Proposition 4.7]{Lack}, an algebra for $Q$ in $\C$ consists of an object $M$ with an $\ifas$-algebra structure and an $\ifas^{op}$-algebra structure subject to compatibility conditions. An $\ifas$-algebra structure is an involutive monoid structure and an $\ifas^{op}$-algebra structure is an involutive comonoid structure. Let $F$ be the $\ifas$-algebra and let $G$ be the $\ifas^{op}$-algebra. The compatibility condition requires that for every bimorphism
\begin{center}
\begin{tikzcd}
\ul{n}\arrow[d, "\phi_1 ", swap]\arrow[r, "f_1"] & \ul{p}\arrow[d, "\phi "]\\
\ul{m}\arrow[r, "f", swap] & \ul{q}
\end{tikzcd}
\end{center}
in the double category $\ifas_2$, the diagram
\begin{center}
\begin{tikzcd}
M^{\otimes n} \arrow[r, "F(f_1)"] & M^{\otimes p}\\
M^{\otimes m} \arrow[u, "G(\phi_1) "] \arrow[r, "F(f)", swap] & M^{\otimes q}\arrow[u, "G(\phi)" , swap]
\end{tikzcd}
\end{center}
commutes. Arguing analogously to \cite[5.3]{Lack}, it suffices to have commutativity for the fundamental bimorphisms of Definition \ref{fund-bimorph}. These are precisely the conditions requiring $M$ to be an involutive bimonoid. Finally we observe that a morphism in $\C$ is a morphism of involutive bimonoids if and only if it preserves the $\ifas$-algebra structure and the $\ifas^{op}$-algebra structure. By \cite[4.8]{Lack} this is true if and only if it preserves the $Q$-algebra structure.
\endproof 

\remark
The theorem tells us that the PROP governing the structure of involutive bimonoids is a composite of the PROPs for involutive monoids and involutive comonoids where the compatibility conditions are precisely the fundamental bimorphisms of the category $\ifas_2$ and the distributive law is determined by the star condition. We have chosen this method of proof as we believe that the double category $\ifas_2$ most neatly encapsulates the structure required to construct the PROP $Q$.

An alternative method of proof would be to define a distributive law between the PROP given in \cite[5.9]{Lack} and the PRO $\mathbb{C}_2$ and analyse the resulting composite. The technical details of such a proof are quite similar to those we have used.
\endremark

\theorem
Let $\C$ be a symmetric monoidal category.
\begin{enumerate}
\item The category $\mathbf{Alg}\left(Q_{\mathcal{V}} , \C\right)$ is equivalent to the category of involutive commutative bimonoids in $\C$.
\item The category $\mathbf{Alg}\left(Q_{\mathcal{H}} , \C\right)$ is equivalent to the category of involutive cocommutative bimonoids in $\C$.
\item The category $\mathbf{Alg}\left(Q_{\mathcal{IF}} , \C\right)$ is equivalent to the category of involutive, commutative, cocommutative bimonoids in $\C$.
\end{enumerate}
\endtheorem
\proof
These equivalences are proved similarly to Theorem \ref{bimon-thm}.
\endproof

\end{document}